\documentclass[11pt,a4paper,oneside]{amsart}
\usepackage{fullpage}
\usepackage{amssymb,amsfonts}
\usepackage{epsfig}
\usepackage{graphicx}

\usepackage{graphics}
\newtheorem{theorem}{Theorem}[section]

\theoremstyle{definition}

\numberwithin{equation}{section}

\DeclareMathOperator{\diam}{diam} 

\newcommand{\be}{\begin{equation}}
\newcommand{\ee}{\end{equation}}

\makeindex

\def\XXint#1#2#3{{\setbox0=\hbox{$#1{#2#3}{\int}$}
 \vcenter{\hbox{$#2#3$}}\kern-.5\wd0}}

\usepackage[colorinlistoftodos]{todonotes}
\usepackage[colorlinks=true, allcolors=blue]{hyperref}

\title{Sobolev embedding implies regularity of measure in metric measure spaces}
\author{Nijjwal Karak}
\address{Department of Mathematical Analysis, Charles University, Sokolovsk\'a 83, 18600 Prague 8, Czech Republic}
\email{nijjwal@gmail.com}
\thanks{This work was supported by OP RDE project no. CZ.02.2.69/0.0/0.0/16\_027/0008495, International Mobility of Researcher at Charles University.}
\begin{document}
\begin{abstract}
We prove that if the Sobolev embedding $M^{1,p}(X)\hookrightarrow L^q(X)$ holds for some $q>p\geq 1$ in a metric measure space $(X,d,\mu),$ then a constant $C$ exists such that $\mu(B(x,r))\geq Cr^n$ for all $x\in X$ and all $0<r\leq 1,$ where $\frac{1}{p}-\frac{1}{q}=\frac{1}{n}.$ This was proved in \cite{Gor17} assuming a doubling condition on the measure $\mu.$ 
\end{abstract}
\maketitle
\indent Keywords: Haj\l asz-Sobolev spaces, Metric measure spaces, Sobolev embeddings.\\
\indent 2010 Mathematics Subject Classification: 46E35.
\section{Introduction}
It is well-known that, for an an open set $\Omega\subset\mathbb{R}^n$ and $1\leq p<n,$ the Sobolev embedding $W^{1,p}(\Omega)\hookrightarrow L^{p^*}(\Omega)$ holds, where $p^*=\frac{np}{n-p},$ if the boundary of $\Omega$ is sufficiently regular, see e.g. \cite{Ada75}. On the other hand, it has been proved in \cite{HKT08} that if the embedding $W^{1,p}(\Omega)\hookrightarrow L^{p^*}(\Omega)$ holds, then $\Omega$ satisfies the so-called measure density condition, i.e. there exists a constant $c>0$ such that for all $x\in X$ and all $0<r\leq 1$
\begin{equation}\label{measuredensity}
\vert B(x,r)\cap\Omega\vert\geq cr^n.
\end{equation}

Let $(X,d,\mu)$ be a metric measure space equipped with a metric $d$ and a Borel regular measure $\mu.$ We assume throughout the note that the measure of every open nonempty set is positive and that the measure of every bounded set is finite. In a metric measure space $(X,d,\mu),$ Haj\l asz \cite{Haj96} has shown that if the space $X$ is $n$-regular, then the embedding $M^{1,p}(X)\hookrightarrow L^{p^*}(X)$ holds, where $p^*=\frac{np}{n-p}.$ Recall that a space $(X,d,\mu)$ is $n$-regular if there exists a constant $C$ such that
\begin{equation}\label{lowerbound}
\mu(B(x,r))\geq Cr^n
\end{equation}
for all $B(x,r)\subset X$ with $r<\diam X.$ Also recall that a $p$-integrable function $u$ belongs to the Haj\l asz-Sobolev space $M^{1,p}(X)$ if there exists a non-negative $g\in L^p(X),$ called a generalized gradient, such that
\begin{equation*}
\vert u(x)-u(y)\vert\leq d(x,y)(g(x)+g(y))\quad\text{a.e. for}\quad x,y\in X.
\end{equation*}
The space $M^{1,p}(X)$ is a Banach space with the norm
\begin{equation*}
\Vert u\Vert_{M^{1,p}(X)}=\Vert u\Vert_{L^p(X)}+\inf \Vert g\Vert_{L^p(X)},
\end{equation*}
where the infimum is taken over all the generalized gradients.\\

In \cite{Gor17}, it has been proved that if the embedding $M^{1,p}(X)\hookrightarrow L^{q}(X)$ holds for some $q>p,$ then the measure $\mu$ satisfies \eqref{lowerbound} for all $x\in X$ and all $0<r\leq 1$ provided that the space $(X,d,\mu)$ is doubling, i.e. there exists a constant $c_d$ such that for every ball $B(x,r),$
$$\mu(B(x,2r))\leq c_d\mu(B(x,r)).$$ 
In this note, we prove the same result but without assuming the doubling condition, as conjectured in \cite{Gor17} and the proof of the same is inspired by \cite{AH} and \cite{KMR15}.  
\begin{theorem}\label{maintheorem}
Let $(X,d,\mu)$ be a metric measure space and $p\geq 1.$ If $M^{1,p}(X)\hookrightarrow L^q(X),$ $q>p,$ then there exists $C=C(p,q,C_e)$ such that
\begin{equation*}
\mu(B(x,r))\geq Cr^n,\quad\text{for}\quad r\in (0,1],
\end{equation*}
where $\frac{1}{p}-\frac{1}{q}=\frac{1}{n}$ and $C_e$ is the constant of the embedding.
\end{theorem}
\section{Proof of Theorem 1.1}
\noindent For each $u\in M^{1,p}(X)$ and for any generalized gradient $g$ of $u$ we have, by the Sobolev embedding,
\begin{equation}\label{embedding}
\left(\int_X\vert u\vert^q\,d\mu\right)^{\frac{1}{q}}\leq C_e\left[\left(\int_X\vert u\vert^p\,d\mu\right)^{\frac{1}{p}}+\left(\int_X g^p\,d\mu\right)^{\frac{1}{p}}\right].
\end{equation}
Fix $x\in X$ and $r\in (0,1].$ For each fixed $j\in\mathbb{N},$ set $r_j=(2^{-j-1}+2^{-1})r,$ and $B_j=B(x,r_j).$ Note that, for all $j\in\mathbb{N},$
\begin{equation*}
\frac{r}{2}<r_{j+1}<r_j\leq\frac{3r}{4}.
\end{equation*}
For each $j\in\mathbb{N},$ let us define $u_j:X\rightarrow\mathbb{R}$ as follows:
\begin{equation*}
u_j(y)=
\begin{cases}
  1& \text{if $y\in B_{j+1}$},\\
  \frac{r_j-d(x,y)}{r_j-r_{j+1}} & \text{if $y\in B_j\setminus B_{j+1}$},\\
  0& \text{if $y\in X\setminus B_j$}.
 \end{cases}
\end{equation*}
It is easy to see that, for each $j\in\mathbb{N},$ $u_j$ is a $(r_j-r_{j+1})^{-1}$-Lipschitz function on $X$ and the function $g_j:=(r_j-r_{j+1})^{-1}\chi_{B_j}$ is a generalized gradient of $u_j.$ In particular, $u_j\in M^{1,p}(X)$ and hence the functions $u_j$ and $g_j$ satisfy \eqref{embedding}. Noting that $(r_j-r_{j+1})^{-1}=2^{j+2}r^{-1}$ we have, for each $j\in\mathbb{N},$
\begin{eqnarray*}
\int_X g_j^p\,d\mu=\frac{2^{p(j+2)}}{r^p}\mu(B_j)\quad
\text{and}\quad
\int_X\vert u_j\vert^p\,d\mu\leq \mu(B_j).
\end{eqnarray*}
Moreover, for each $j\in\mathbb{N},$
\begin{equation*}
\int_X\vert u_j\vert^q\,d\mu\geq \mu(B_{j+1}).
\end{equation*}
Use these estimates while applying \eqref{embedding} for the pair $(u_j,g_j),$ for every $j\in\mathbb{N},$ to obtain
\begin{eqnarray}\label{fromembedding}
\mu(B_{j+1})^{1/q} &\leq & C_e\left(1+\frac{2^{j+2}}{r}\right)\mu(B_j)^{1/p} \nonumber \\
&\leq & \frac{C_e}{r}2^{j+3}\mu(B_j)^{1/p},
\end{eqnarray}
where in the last inequality we have used the fact that $r\leq 1.$ Raising both sides of the inequality \eqref{fromembedding} to the power $p/\alpha^{j-1},$ where $\alpha=q/p\in (1,\infty),$ yields
\begin{equation*}
\mu(B_{j+1})^{1/\alpha^j}\leq\left(\frac{C_e}{r}\right)^{p/\alpha^{j-1}}2^{p(j+3)/\alpha^{j-1}}\mu(B_j)^{1/\alpha^{j-1}}.
\end{equation*}
Letting $P_j=\mu(B_j)^{1/\alpha^{j-1}},$ we rewrite the above inequality as
\begin{equation*}
P_{j+1}\leq \left(\frac{C_e}{r}\right)^{p/\alpha^{j-1}}2^{p(j+3)/\alpha^{j-1}}P_j\quad\forall j\in\mathbb{N}.
\end{equation*}
After iteration, we obtain, for every $j\in\mathbb{N},$
\begin{equation}\label{product}
P_{j+1}\leq P_1\prod_{k=1}^j2^{p(k+3)/\alpha^{k-1}}\left(\frac{C_e}{r}\right)^{p/\alpha^{k-1}}.
\end{equation}
Observe that
\begin{equation*}
\prod_{k=1}^{\infty}\left(\frac{C_e}{r}\right)^{p/\alpha^{k-1}}=\left(\frac{C_e}{r}\right)^{p\sum_{k=1}^{\infty}\alpha^{1-k}}=\left(\frac{C_e}{r}\right)^{\frac{p\alpha}{\alpha-1}}
\end{equation*}
and
\begin{equation*}
\prod_{k=1}^{\infty}2^{p(k+3)/\alpha^{k-1}}=2^{p\sum_{k=1}^{\infty}(k+3)\alpha^{1-k}}=2^{\frac{p\alpha^2}{(\alpha-1)^2}+\frac{3p\alpha}{\alpha-1}}.
\end{equation*}
On the other hand, from the construction of $B_j$'s, we have
$$\mu(B(x,r/2))^{1/\alpha^{j-1}}\leq P_j=\mu(B_j)^{1/\alpha^{j-1}}\leq \mu(B(x,r))^{1/\alpha^{j-1}}$$
and therefore $\lim_{j\rightarrow\infty}P_j=1.$ Consequently, passing to the limit in \eqref{product} and using $P_1\leq \mu(B(x,r)),$ we obtain
\begin{equation*}
1\leq 2^{\frac{p\alpha^2}{(\alpha-1)^2}+\frac{3p\alpha}{\alpha-1}}\left(\frac{C_e}{r}\right)^{\frac{p\alpha}{\alpha-1}}\mu(B(x,r)).
\end{equation*}
Therefore
\begin{equation*}
\mu(B(x,r))\geq Cr^{\frac{p\alpha}{\alpha-1}},
\end{equation*}
where
\begin{equation*}
\frac{1}{C}=2^{\frac{p\alpha^2}{(\alpha-1)^2}+\frac{3p\alpha}{\alpha-1}}C_e^{\frac{p\alpha}{\alpha-1}}.
\end{equation*}
Finally, we use $q=np/(n-p)$ and $\alpha=q/p$ to get the desired result.\qedhere
\def\bibname{References}
\bibliography{embedding_measure}

\begin{thebibliography}{1}

\bibitem{Ada75}
Robert~A. Adams.
\newblock {\em Sobolev spaces}.
\newblock Academic Press [A subsidiary of Harcourt Brace Jovanovich,
  Publishers], New York-London, 1975.
\newblock Pure and Applied Mathematics, Vol. 65.

\bibitem{AH}
Ryan Alvarado and Piotr Haj{\l}asz.
\newblock A note on metric-measure spaces supporting {P}oincar{\'e}
  inequalities.
\newblock preprint 2019, https://arxiv.org/abs/1902.10876.

\bibitem{Gor17}
Przemys{\l}aw G{\'o}rka.
\newblock In metric-measure spaces {S}obolev embedding is equivalent to a lower
  bound for the measure.
\newblock {\em Potential Anal.}, 47(1):13--19, 2017.

\bibitem{Haj96}
Piotr Haj{\l}asz.
\newblock Sobolev spaces on an arbitrary metric space.
\newblock {\em Potential Anal.}, 5(4):403--415, 1996.

\bibitem{HKT08}
Piotr Haj{\l}asz, Pekka Koskela, and Heli Tuominen.
\newblock Sobolev embeddings, extensions and measure density condition.
\newblock {\em J. Funct. Anal.}, 254(5):1217--1234, 2008.

\bibitem{KMR15}
Lyudmila Korobenko, Diego Maldonado, and Cristian Rios.
\newblock From {S}obolev inequality to doubling.
\newblock {\em Proc. Amer. Math. Soc.}, 143(9):4017--4028, 2015.

\end{thebibliography}
\bibliographystyle{plain}
\end{document}